\documentclass[a4paper,12pt]{amsart}
\usepackage{amsfonts,amssymb,epsfig}
\usepackage[latin1]{inputenc}
\newtheorem{thm}{Theorem}[section]
\newtheorem{cor}[thm]{Corollary}
\newtheorem{lem}[thm]{Lemma}
\newtheorem{pro}[thm]{Proposition}
\newtheorem{dfn}[thm]{Definition}
\newtheorem{rmq}[thm]{Remark}
\newtheorem{expl}[thm]{Example}

\def\dessous#1\sous#2{\mathrel{\mathop{\kern0pt#2}\limits_{#1}}}
\newcommand{\mbb}[1]{\mathbb{#1}}

\newcommand{\nb}{\nonumber}

\newcommand{\noi}{\noindent}
\newcommand{\n}{\nabla}

\newcommand{\beq}{\begin{eqnarray}}
\newcommand{\eeq}{\end{eqnarray}}
\newcommand{\beqs}{\begin{eqnarray*}}
\newcommand{\eeqs}{\end{eqnarray*}}

\newcommand{\bpro}{\begin{pro}}
\newcommand{\epro}{\end{pro}}
\newcommand{\blem}{\begin{lem}}
\newcommand{\elem}{\end{lem}}
\newcommand{\bdfn}{\begin{dfn}}
\newcommand{\edfn}{\end{dfn}}
\newcommand{\bcor}{\begin{cor}}
\newcommand{\ecor}{\end{cor}}
\newcommand{\bthm}{\begin{thm}}
\newcommand{\ethm}{\end{thm}}
\newcommand{\bex}{\begin{expl}}
\newcommand{\eex}{\end{expl}}
\newcommand{\brmq}{\begin{rmq}}
\newcommand{\ermq}{\end{rmq}}
\begin{document}
\title{Biharmonic Reeb curves in   Sasakian manifolds}
\author{S. Degla}
\address{Institut de Mathematiques et de Sciences Physiques, 01 BP 613, Porto-Novo, Benin.}
\email{sdegla@imsp-uac.org}
\author{L. Todjihounde}
\email{leonardt@imsp-uac.org}

\date{}


\begin{abstract}
Sasakian manifolds provide explicit formulae of some Jacobi operators which describe the biharmonic equation of curves in 
Riemannian manifolds. In this paper we characterize non-geodesic biharmonic curves in Sasakian manifolds which are either tangent
 or normal to the Reeb vector field. \\In the three-dimensional case, we prove that such curves are some helixes whose geodesic curvature 
and geodesic torsion satisfy a given relation.
\end{abstract}

\maketitle
\section{Introduction}
The notions of harmonic and biharmonic maps between Riemannian manifolds have been introduced by
J. Eells and J.H. Sampson  (see \cite{ES}). \\ For a map 
 $\phi:\, (M,g)\rightarrow (N,h)$ between Riemannian manifolds
the energy functional $E_1$ is defined by
$$E_1(\phi) \ =\ \frac{1}{2}\int_M |d\phi|^2v_g.$$

\noindent Critical points of $E_1$ are called harmonic maps and are then solutions of the corresponding
Euler-Lagrange equation  $$\tau_1(\phi)=trace\nabla^\phi 
 d\phi .$$ 
Here $\nabla^\phi$ denotes the induced connection on the pull-back
 bundle $\phi^{-1}(TN)$ and $\tau_1(\phi)$ is called the tension field of $\phi$.

\noindent Biharmonic maps are the critical points of the functional bienergy 
$$
E_2(\phi) \ =\ \frac{1}{2}\int_M |\tau_1(\phi)|^2v_g,
$$
whose Euler-Lagrange equation is given by the vanishing of the bitension field
(cf. \cite{Y2}) defined by
 \beq\label{o1}
 \tau_{2}(\phi)
               & =&-\Delta^\phi \tau_{1}(\phi)
                    -trace R^N(d\phi, \tau_{1}(\phi))d\phi,
\eeq
where $\Delta^\phi=-trace_g(\n^\phi\n^\phi-\n^\phi_\n)$ is the Laplacian on the sections of $\phi^{-1}(TN)$, and $R^N$ is
the Riemannian curvature operator of $(N, h)$. Note that 
\beq \tau_2(\phi)=J_\phi(\tau_1(\phi))\eeq
where $J_\phi$ is the Jacobi operator along $\phi$ defined by 
\beq
J_\phi(X)=-\Delta^\phi X
                    -trace R^N(d\phi, X)d\phi,\quad \forall\ X \ \in\ \phi^{-1}(TN).
\eeq
Harmonic maps are obviously biharmonic and are absolute minimum of the bienergy.

\noindent Nonminimal biharmonic submanifolds of the pseudo-euclidean space and of the spheres have been studied in \cite{CI} and \cite{CMO2}.

\noindent Biharmonic curves have been investgated on many special Riemannian manifolds like Heisenberg groups \cite{COP}, \cite{DF2}, invariant surfaces 
\cite{MO}, Damek-Ricci spaces \cite{DT}, Sasakian manifolds  \cite{DF}, etc.

\noindent As in the general theory of metric contact manifolds an important role is played by the Reeb vector field whose dynamics can be used
to study the structure of the contact manifold or even the underlying manifold using techniques of Floer homology such as
symplectic field theory and embedded contact homology.

\noindent The main purpose of this work is to study non-geodesic biharmonic curves in a ($2n+1$)-dimensional Sasakian manifold, which are
 either tangent or normal to the Reeb vector field.
\section{Sasakian manifolds}
\noindent A contact manifold is a $(2n+1)$-dimensional manifold $M$ equipped with a global $1$-form $\eta$ such that 
$\eta \wedge (d\eta)^n\neq 0$ everywhere on $M$. It has an underlying almost contact structure $(\eta,\varphi,\xi)$ where
$\xi$ is a global vector field (called the characteristic vector field) and $\varphi$ a global tensor of type $(1,1)$ such that
\beq
\eta(\xi)=1,\ \varphi\xi=0,\ \eta\varphi=0,\ \varphi^2= -I+\eta\otimes\xi.
\eeq
A Riemannian metric $g$ can be found such that
\beq
\eta=g(\xi, .), \ d\eta=g(.,\varphi .),\ g(.,\varphi .)=-g(\varphi ., .).
\eeq
$(M,\eta,g)$ or $(M,\eta,g,\xi,\varphi)$ is called a contact metric  manifold. If the almost complex structure $J$
on $M\times \mbb{R}$ defined by 
\beq
J(X,f\frac{d}{dt})=(\varphi X-f\xi,\eta(X)\frac{d}{dt}),
\eeq
is integrable, $(M,\eta,g)$ is said to be Sasakian.

\noindent The following relations play an important role in the present work:
\blem\label{l1}\cite{DB2}
On a Sasakian manifold $(M,\eta,g)$ we have
\beq
R(X,\xi)X = -\xi
\eeq
and
\beq
R(\xi,X)\xi=-X
\eeq
for any unit vector field $X$ orthogonal to the Reeb vector field $\xi$, where $R$ denotes the Riemannian curvature of $(M,g)$.
\elem
\section{Biharmonic curves in  Sasakian manifolds}
\noi Let
$ \gamma:\ I\ \longrightarrow \ (M^{2n+1},\eta, g)$ be a regular curve
 parametrized by its arc lenght in a (2n+1)-dimensional Sasakian manifold and  $\{T,N_1, ...,N_{2n}\}$ be the Frenet frame in
$M^{2n+1}$, defined along $\gamma$, where $T=\gamma'$ is  the unit tangent  vector field of $\gamma$. \\ It holds:
\blem\cite{DF}\label{l2}
The Frenet equations of $\gamma$ are given by
$$\left\{\begin{array}{lll}\n_TT&=& \chi_1 N_1\\
          \n_TN_1& =& -\chi_1 T +\chi_2 N_2\\
           \vdots &  \vdots&    \vdots \\
        \n_TN_k & =& -\chi_k N_{k-1} +\chi_{k+1}N_{k+1},\quad k=2,...,2n-1,\\
\vdots &  \vdots&    \vdots \\
     \n_TN_{2n} &=& -\chi_{2n} N_{2n-1},
          \end{array}  \right. $$
where $\chi_1=|\n_TT|$, $\chi_2=\chi_2(s),...,\chi_{2n}=\chi_{2n}(s)$ are real valued functions, where $s$ is the arc length
of $\gamma$. 
\elem
\bdfn
If the functions $\chi_k$, $k=1,...,2n$ are all constant, then $\gamma$ is said to be a helix..
\edfn
\noindent The tension field $\tau_1(\gamma)$ and the bitension field $\tau_2(\gamma)$ of the curve $\gamma$ are given in the 
Frenet frame $(T,N_1,...,N_{2n})$ by:
\bpro\label{pr3}\cite{DF}
\beq\label{e1}
\tau_1(\gamma)= \chi_1 N,
\eeq
and
\beq\label{e2}
\tau_2(\gamma)&=&-3\chi_1\chi_1' T+(\chi_1''-\chi_1^3-\chi_1\chi_2^2)N_1\\ 
              &&-(2\chi_1'\chi_2+\chi_1\chi_2')N_2 +\chi_1\chi_2\chi_3N_3-\chi_1 R(T,N_1)T.\nb
\eeq
\epro
\noindent From Proposition \ref{pr3}, we get:
\bpro\label{pr4}
If $\gamma$ is either tangent or normal to the Reeb vector field, then
\beq
\tau_2(\gamma)&=&-3\chi_1\chi_1' T+(\chi_1''-\chi_1^3-\chi_1\chi_2^2+\chi_1)N_1\label{e04}\\ 
              &&-(2\chi_1'\chi_2+\chi_1\chi_2')N_2 +\chi_1\chi_2\chi_3N_3 \nb
\eeq
\epro
Proof\\
Let $\gamma$ be non-geodesic biharmonic curve in a Sasakian manifold $(M,\eta,g)$. \\
Assume that $\gamma$ is tangent to the Reeb vector field $\xi$; that is $T=\xi$. \\ The relation (\ref{e2}) in Proposition
\ref{pr3} becomes
\beq
\tau_2(\gamma)&=&-3\chi_1\chi_1' T+(\chi_1''-\chi_1^3-\chi_1\chi_2^2)N_1\nb \\ \nb
              &&-(2\chi_1'\chi_2+\chi_1\chi_2')N_2 +\chi_1\chi_2\chi_3N_3-\chi_1 R(\xi,N_1)\xi\\ \nb
              &=& -3\chi_1\chi_1' T+(\chi_1''-\chi_1^3-\chi_1\chi_2^2)N_1\\ \nb
              &&-(2\chi_1'\chi_2+\chi_1\chi_2')N_2 +\chi_1\chi_2\chi_3N_3+\chi_1N_1, \mbox{ according to lemma } \ref{l1}\\ \nb
\eeq
It follows that
\beq
\tau_2(\gamma)&=&-3\chi_1\chi_1' T+(\chi_1''-\chi_1^3-\chi_1\chi_2^2+\chi_1)N_1\\ \nb 
              &&-(2\chi_1'\chi_2+\chi_1\chi_2')N_2 +\chi_1\chi_2\chi_3N_3
\eeq
We assume now that $\gamma$ is normal to the Reeb vector field; that is $N_1=\xi$. \\
The relation (\ref{e2}) in Proposition 
\ref{pr3} becomes then
\beq
\tau_2(\gamma)&=&-3\chi_1\chi_1' T+(\chi_1''-\chi_1^3-\chi_1\chi_2^2)N_1 \nb \\ \nb
              &&-(2\chi_1'\chi_2+\chi_1\chi_2')N_2+\chi_1\chi_2\chi_3N_3-\chi_1 R(T,\xi)T\\ \nb
              &=&-3\chi_1\chi_1' T+(\chi_1''-\chi_1^3-\chi_1\chi_2^2)N_1\\ \nb
              &&-(2\chi_1'\chi_2+\chi_1\chi_2')N_2+\chi_1\chi_2\chi_3N_3+\chi_1 \xi \mbox{ according to lemma } \ref{l1}\\ \nb
\eeq
We obtain then again
\beq
\tau_2(\gamma)&=&-3\chi_1\chi_1' T+(\chi_1''-\chi_1^3-\chi_1\chi_2^2+\chi_1)N_1\\ \nb 
              &&-(2\chi_1'\chi_2+\chi_1\chi_2')N_2 +\chi_1\chi_2\chi_3N_3              
\eeq
Thus we get the relation (\ref{e04}) in both cases.

$\hfill{\square}$

\noindent From Proposition \ref{pr4}, we get the following result.
\bthm
Non-geodesic biharmonic curves in Sasakian
manifolds which are either tangent or normal to the Reeb vector field are characterized by  :
\beq\label{e03}
\left\{ 
\begin{array}{l}
 \chi_1 = constant \in [-1,0[U]0,1],\\
\chi_2=\pm\sqrt{1-\chi_1^2},\\
\chi_2\chi_3 = 0.
\end{array}
\right.
\eeq
where $\chi_1$, $\chi_2$ and $\chi_3$ are functions defined in lemma \ref{l2}.
\ethm 
Proof\\
 $\tau_2(\gamma)=0$ with $\chi_1\neq 0$ implies $\chi_1'=0$ according to the first component in (\ref{e04}). So $\chi_1$ is
constant. Then the second component gives $\chi_1^2+\chi_2^2=1$.    Thus (\ref{e03}) is satisfied.

%
                                                                                             $\hfill{\square}$
\bcor
If $\chi_1=\pm 1$ then $\chi_2=0$. And if $\chi_1\neq \pm 1$ then $\chi_3=0$.
\ecor
\brmq
From the conditions given in (\ref{e03}), it is clear that in general  non-geodesic biharmonic curves in Sasakian manifolds
are not helixes since only the functions $\chi_1$ and $\chi_2$ have to be constant and maybe $\chi_3$.
\ermq
\noindent In three-dimensional Sasakian manifolds the Frenet frame is given by 
$$\left\{\begin{array}{lll}\n_TT&=& \chi_1 N_1\\
          \n_TN_1& =& -\chi_1 T +\chi_2 N_2\\
           \n_TN_2 & =& -\chi_2 N_1,
          \end{array}  \right. $$
where $\chi_1=|\n_TT|$, $\chi_2=\chi_2(s)$ is a real valued function, where $s$ is the arc length
of $\gamma$. \\
And the equation characterizing the non-geodesic biharmonic curves which are either tangent or normal to Reeb vector field is 
reduced to
   \beq\label{e5}
\tau_2(\gamma)&=&-3\chi_1\chi_1' T+(\chi_1''-\chi_1^3-\chi_1\chi_2^2+\chi_1)N_1\\ 
              &&-(2\chi_1'\chi_2+\chi_1\chi_2')N_2 .\nb
\eeq  
So we have the following result.

\bthm
Non-geodesic biharmonic curves which are either tangent or normal to the Reeb vector field in  three-dimensional Sasakian
manifolds are helixes whose geodesic curvature $\chi_1$ and geodesic torsion $\chi_2$ are related by: 
$$\chi_1^2+\chi_2^2=1,\ \mbox{with}\ \chi_1\neq 0 \ .$$
\ethm

%

\bibliographystyle{plain}

\end{document}